\input amstex
\documentstyle{amsppt}
\magnification=\magstep1 \baselineskip=18pt \hsize=6truein
\vsize=8truein

\topmatter
\title On a conformal gap and finiteness theorem for a
class of four manifolds
\endtitle
\author Sun-Yung A. Chang$^{\ssize 1}$, Jie Qing$^{\ssize 2}$
and Paul Yang$^{\ssize 3}$
\endauthor

\adjustfootnotemark{1}\footnotetext{Princeton University, Department
of Mathematics, Princeton, NJ 08544-1000, supported by NSF Grant
DMS--0245266.} \adjustfootnotemark{1}\footnotetext{University of
California, Department of Mathematics, Santa Cruz, CA 95064,
supported in part by NSF Grant DMS--0402294.}
\adjustfootnotemark{1}\footnotetext{Princeton University, Department
of Mathematics, Princeton, NJ  08544-1000, supported by NSF Grant
DMS--0245266.}

\vskip .1in

\leftheadtext{gap and finiteness theorems } \rightheadtext{Chang,
Qing and Yang}

\abstract In this paper we develop a bubble tree structure for a
degenerating class of Riemannian metrics satisfying some global
conformal bounds on compact manifolds of dimension 4.  Applying the
bubble tree structure, we establish a gap theorem, a finiteness
theorem for diffeomorphism type for this class, and make a
comparison of the solutions of the $\sigma _k$ equations on a
degenerating family of Bach flat metrics.

\endabstract

\endtopmatter

\def\bcw{\mathbin{\bigcirc\mkern-15mu\wedge}}

\document
\noindent{\bf 1. Introduction} \vskip 0.1in

Given a Riemannian four-manifold $(M^4,g)$, let $Rm$ denote the
curvature tensor, $W$ the Weyl curvature tensor, $Ric$ the Ricci
tensor, and $R$ the scalar curvature.  The usual decomposition of
$Rm$ under the action of $O(4)$ can be written
$$
Rm = W + {\tsize{1\over 2}} E  \bcw g + {\tsize{1\over 24}} R g \bcw
g, \tag1.1
$$
where $E = Ric - {1 \over 4} Rg$ is the trace-free Ricci tensor and
$\bcw$ denotes the Kulkarni-Nomizu product. We also recall the
Gauss-Bonnet-Chern formula:
$$
8 \pi^2 \chi(M^4) = \int_{M^4} \bigl( {\tsize{1\over 4}}|W|^2 +
    {\tsize{1\over 24}}R^2  \, - \,
    {\tsize{1\over 2}} |E|^2 \bigr)  dvol.
\tag1.2
$$
In [CGY-2], based on an earlier work of [Ma], we have the following
sharp form of "sphere" theorem:

\proclaim{Theorem 1.1} Let $(M^4,g)$ be a smooth, closed manifold
for which \roster
\item"{(i)}" the Yamabe constant
$Y(M^4, [g]) > 0$, and
\smallskip
\item"{(ii)}" the curvature satisfies
$$
\int_{M^4} \bigl(
       {\tsize{1\over 24}}R^2 -{\tsize{1\over2}} |E|^2 \bigr) dvol  > \int_{M^4} {\tsize{1\over4}} |W|^2  dvol \geq  0.
\tag1.3
$$
\endroster
\smallskip\noindent
Then $M^4$ is diffeomorphic to either $S^4$ or ${\Bbb {RP}}^4$.
Furthermore, if $M$ is not diffeomorphic to $(S^4, g_c)$ or $({\Bbb
{RP}}^4, g_c)$ and the inequality in (1.3) becomes an inequality;
then either \newline \noindent (1) $(M, g)$ is conformally
equivalent to $({\Bbb {CP}}^2, g_{FS})$, or
\newline \noindent  (2) $ (M, g)$ is conformal equivalent to
$( (S^3 \times S^1)/\Gamma, g_{prod})$ for a finite group $\Gamma$.
\endproclaim

\remark{Remarks} \roster
\item"{1.}"  Recall that the Yamabe constant is defined by
$$
Y(M^4, [g]) = \inf_{\tilde g \in [g]}
    vol(\tilde g)^{-{1\over 2}}\int_{M^4} R_{\tilde g} dvol_{\tilde g},
$$
where $[g]$ denotes the conformal class of $g$. Positivity of the
Yamabe invariant implies that $g$ is conformal to a metric of
strictly positive scalar curvature.

\item"{2.}"  In the statement of Theorem 1.1 the norm of the Weyl
tensor is given by $|W|^2 = W_{ijkl}W^{ijkl}$;
\endroster
\endremark

The proof of the second part of the Theorem above relies on the
vanishing of the Bach tensor. If $g$ satisfies the condition of
Theorem 1.1 and $(M^4, g)$ is not diffeomorphic to either  $S^4$ or
${\Bbb {RP}}^4$, then
 $g$ is a critical point (actually, a local minimum)
of the {\it Weyl functional} $g \mapsto \int |W|^2 dvol$. The
gradient of this functional is called the {\it Bach tensor}; which
we shall define in Section 2 below, and we will say that critical
metrics are {\it Bach flat}.  Note that the conformal invariance of
the Weyl functional in dimension four implies that Bach-flatness is
a conformally invariant property.

We now define some notations to state results in this paper. Given a
Riemannian four--manifold $(M^4, g)$, the {\it Weyl--Schouten}
tensor is defined by
$$
A = Ric - {1\over 6} Rg
$$
In terms of the Weyl--Schouten tensor, the decomposition $(1.1)$ can
be written as
$$
\text{Rm} = W + {1\over 2} A \bcw g \,. \tag1.4
$$
This splitting of the curvature tensor induces a splitting of the
Euler form. To describe this, we introduce the elementary symmetric
polynomials
$$
\sigma_\kappa ( \lambda_1 , \ldots , \lambda_n) =
   \sum\limits_{i_1 < \cdots < i_k}
     \lambda_{i_1} \cdots \lambda_{i_\kappa} \, ;
$$
where ${\lambda_i}'s$ denote the eigenvalues of the (contracted)
tensor $g^{-1} A$. To simplify notations, we denote $\sigma_\kappa
(A_g) = \sigma_\kappa (g^{-1}A)$. We note for a manifold of
dimension 4,
$$
\sigma_2 (A_g) = {\tsize{1\over 24}}R^2 -{\tsize{1\over2}} |E|^2.
\tag 1.5
$$
Sometimes we will also denote $\sigma_2 (A_g)$ as $\sigma_2$ when
the metric $g$ is fixed.

The Gauss--Bonnet--Chern formula (1.2) may be written as
$$
8 \pi ^2 \chi (M^4) =
   \int_{M^4} {1\over 4 }|W|^2 dvol +
     \int_{M^4} \sigma_2 dvol \;.
\tag 1.6
$$
Note that the conformal invariance of the Weyl functional implies
that the quantity
$$
\int_{M^4} \sigma_2 dvol
$$
is conformally invariant as well.

For a four manifold with a positive Yamabe constant, it follows from
the solution of the Yamabe problem ([Au], [S]) that we may assume
that $g$ is the Yamabe metric which attains the Yamabe constant,
then $R_g$ is a constant and
$$
\aligned \int_M \sigma_2 (A_g) dv_g & \leq \int_M \frac 1{24}R_g^2
dv_g \,
 = \,  \frac 1{24} \frac {(\int_M R_g dv_g)^2}{ vol(g) } \\
& \leq \frac 1{24} \frac {(\int_M R_c dv_{g_c})^2}{ (vol(g_c)) } \, = \, 16 \pi^2; \\
\endaligned
\tag 1.7
$$
and equality holds if and only if $(M^4, g)$ is conformally
equivalent to the standard 4-sphere $(S^4, g_c)$ with $R_c = R_{g_c}
= 12$ and $ vol(g_c) = \frac { 8 \pi^2}{3} $.

In view of the inequality (1.7), it is natural to ask whether a
stronger form of rigidity holds in the statement of Theorem 1.1 in
which the integral of $\sigma_2$ is compared directly to the
constant $ 16 \pi^2$ instead to the $L^2$ integral of the Weyl
tensor. Such a possibility is suggested by the recent remarkable
rigidity result of Bray-Neves [BN] on compact manifolds of dimension
3. To state their result, recall the Yamabe invariant is defined as
$Y(M) = sup_{g} Y(M, [g])$, and denote the Yamabe constants $Y_1 =
Y(S^3, [g_c])$, $Y_2 = Y(RP^3, [g_c])$.

\proclaim{Theorem 1.2 ([BN])} A closed 3-manifold with Yamabe
invariant $Y > Y_2$ is either $S^3$ or a
 connect
sum with a $S^2$ bundle over $S^1$.
\endproclaim

The gap theorem in this paper is a first step in this direction:

\proclaim{Theorem A} Suppose that $(M^4, g)$ is a Bach flat closed
4-manifold with positive Yamabe constant and that
$$
\int_M |W|_g^2dv_g  \leq \Lambda_0. \tag 1.8
$$
for some fixed positive number $\Lambda_0$. Then there is a positive
constant $\epsilon_0$ such that if
$$
\int_M \sigma_2 (A_g) dv_g \geq (1-\epsilon) 16\pi^2 \tag 1.9
$$
holds for some constant $\epsilon < \epsilon_0$, then $(M^4, g)$ is
conformally equivalent to the standard 4-sphere.
\endproclaim

The analysis that is developed to prove the gap theorem above can be
adopted to prove a result of finiteness of diffeomorphism classes of
manifolds satisfying some conformally invariant conditions and $L^2$
bounds:

\proclaim{Theorem B} Suppose that $\Cal A$ is a collection of Bach
flat Riemannian manifolds $(M^4, g)$ with positive Yamabe constants,
satisfying inequality (1.8) for some fixed positive number
$\Lambda_0$, and that
$$
\int_M (\sigma_2 (A_g) dv_g \geq a_0, \tag 1.10
$$
for some fixed positive number $a_0$. Then there are only finitely
many diffeomorphism types among manifolds in $\Cal A$.
\endproclaim

It is known that in each conformal class of metrics belonging to the
family $\Cal A$, there is a metric $\bar g=e^{2w}g$ such that
$\sigma _2(A_{\bar g}) =1$, which we shall call the $\sigma _2$
metric. The recent work of Gursky and Viaclovsky [GV] also showed
that if in addition, the positive cone $\Gamma _k^+$ is nonempty for
$k=3$ or $4$, then there exists a conformal metric with $\sigma
_k(A_{\bar g}) =1$. We shall call these metrics $\sigma _k$ metrics.
As an application of the bubbling tree analysis, we will show:

\proclaim{Theorem C} For the conformal classes $\left[g_0\right] \in
{\Cal A}$ the conformal metrics $g=e^{2w}g_0$ satisfying the
equation $\sigma _2(g)=1$ has a uniform bound for the diameter.
\endproclaim

In cases where the $\sigma _k=1$ metrics exists for $k>2$, the Ricci
tensor has a positive apriori lower bound, hence the diameter bound
follows immediately.

\noindent {\bf Remarks:} \vskip .1in

\noindent 1. Our proof of the theorems above builds upon some
estimates in the recent work of Tian-Viaclovsky ([TV-1], [TV-2]; see
also [An-2]) on the compactness of Bach-flat metrics on 4-manifolds;
but our proof relies on a finer analysis of the concentration
phenomenon near points of curvature concentrations. To do this, we
build a bubble tree consisting of vertices which are bubbles near
points of concentrations, and edges consisting of neck regions
connecting different vertexes (see section 4 for a more precise
definition of the bubble trees.) The idea of using bubble tree
construction to achieve finite diffeomorphism types for classes of
manifolds under suitable curvature conditions was developed in the
work of Anderson-Cheeger [AC]. Our construction is modeled after
this work but differs in the way that our bubble tree is built from
the bubbles at points with the smallest scale of concentration, that
is around those points {p} with the smallest radius $r$ so that the
geodesic ball $B_r(p)$ centered around $p$ with radius $r$ achieved
some fixed energy $\int_{B_r(p)}|Rm|^2 dvol $, to bubbles with
larger scale; while the tree in [AC] is constructed from bubbles of
large scale to bubbles with smaller scales. The inductive method of
construction of our bubble tree is modeled on earlier work of [BC],
[Q] and [St] on the study of concentrations of energies in harmonic
maps and the scalar curvature equations.

\vskip .1in \noindent 2. Our proof does not use the stronger volume
estimates, that is, the uniform volume growth for any geodesic ball
in a Bach flat 4-manifold with positive Yamabe constant, $L^2$
bounds of curvature, and bounded first Betti number developed in
[TV-2] and [An-2].  Instead we obtain as a corollary (see Corollary
4.6 in section 4 below) of our bubble tree construction some uniform
estimates for the {\it intrinsic} diameter of the geodesic spheres
near points of curvature concentrations for this class of manifolds.
Indeed one can derive the uniform volume growth as a consequence of
the uniform bound of the diameters. \vskip .1in \noindent 3. Since
the neck theorem and the bubble tree construction which we have
derived in this paper are not "uniform" in scale, we cannot derive
our version of the result of finite diffeomorphism type (Theorem B)
by directly applying the arguments in [AC]. Instead we have
established the proof of Theorem B by an argument of contradiction.
\vskip .1in \noindent 4. Since smooth points of the limit space may
possibly be points of curvature concentration, the proof of Lemma
2.16 in [AC] is not valid in our setting. \vskip .1in\noindent

This paper is organized as follows. In section 2, we discuss some
preliminaries and recall some results in [CGY-2] and [G], and some
key estimates in [TV-1].  In section 3, we present a neck theorem,
which is a variant form of the neck theorem in [AC]. In section 4 we
describe our bubble tree construction; which is the major part of
the paper. In section 5 we apply the bubble tree construction to
prove our main results Theorem A and Theorem B. In section 6, we
develop the analysis of the $\sigma _2$ equation on an ALE bubble.

In a subsequent paper, we will apply the result in Theroem C of this
paper to further study some classification problem of metrics in
$\Cal A$.

The authors wish to thank Tian and Viaclovsky for informative
discussions of the subject.

\vskip 0.1in \noindent{\bf 2. Preliminaries}\vskip 0.1in

In this section we recall some known facts and quote some recent
work in [TV-1] which forms the basis for our construction of the
bubble tree in section 4.

Recall our setting is a
 compact, closed four manifold
$(M^4, g)$ with positive Yamabe constant $Y(M, [g])$. One basic fact
([Au], [S]) is that on such a manifold, $Y(M, [g])$ is achieved by a
metric,
 called the Yamabe metric and which we
 denote again by $g$, with constant scalar curvature
$R_g$.  We also have the following Sobolev inequality:

\proclaim{Lemma 2.1} \quad Suppose that $(M^4, g)$ is a closed,
compact 4-manifold with positive Yamabe constant and suppose that
$g$ is a Yamabe metric. Then
$$
Y(M, [g]) (\int_M u^4 dv_g)^\frac 12 \leq 6 \int_M |\nabla_g u|^2
dv_g + \int_M R_g u^2 dv_g, \tag 2.1
$$
for all function $u \in W^{1,2}(M)$.
\endproclaim

Notice that the inequality (2.1) is invariant under the scaling of
metrics $ g \rightarrow cg$. Thus we may consider 4-manifolds which
are complete, non-compact and have infinite volume and are limits of
rescaled manifolds which satisfy conditions in Lemma 2.1 with a
common lower bound on the Yamabe constants. On such a manifold, the
Sobolev inequality (2.1) takes the following form:
$$
(\int_M u^4 dv_g)^\frac 12 \leq  C_s \int_M |\nabla_g u|^2 dv_g \tag
2.2
$$
for any $u \in C^1(M)$ with compact support inside $M$.

As a by-product of the Sobolev inequality (2.2), one can derive a
lower bound on the volume of all geodesic balls defined on $M$. (cf.
Lemma 3.2 in [He].)

\proclaim{Lemma 2.2} Suppose that either the Sobolev inequality
(2.1) or (2.2) hold on $(M^4, g)$, then there exists a positive
number $v_0$ depending only on $C_s$ such that
$$
\text{vol} (B_r(p)) \geq v_0 r^4 \tag 2.3
$$
for all $p\in M$ and $r < \text{dist}(p, \partial M)$, where
$B_r(p)$ is the geodesic ball of radius $r$ centered at $p$.
\endproclaim

Recall (cf [De]) in local coordinate, the Bach tensor is defined as
$$
B_{ij} \, = \, \nabla^k \nabla^\ell W_{kij\ell} \, + \, \frac{1}{2}
\, R^{k \ell} \, W_{k i j \ell}. \tag 2.4
$$
Using the Bianchi identities, this can be rewritten as
$$
\align B_{ij} \, = \, & - \frac{1}{2} \Delta E_{ij} \, +
\frac{1}{6} \, \nabla_i \nabla_j R \, - \frac{1}{24} \, \Delta R
g_{ij} - E^{k \ell}
W_{ikj\ell} \, \, \\
& + E_i^k E_{jk} \, - \frac{1}{4} \, |E|^2 g_{ij} \, + \,
\frac{1}{6} \,
R E_{ij}, \tag 2.5  \\
\endalign
$$
where $\Delta E_{ij} = g^{kl} \nabla_{k} \nabla_{l} E_{ij}$.

A Bach flat metric is a metric in which the Bach tensor
vanishes--this happens for example when the metric a critical metric
for the functional which is the $L^2$ norm of the Weyl tensor. The
following are some basic identities for Bach flat metrics on
manifolds of dimension 4.

\proclaim{Lemma 2.3} (See [De], and Proposition 3.3 in [CGY-2]). If
$(M^4 , g)$ is Bach flat, then
$$
0 = \int\limits_{M^4} \left\{ 3 \, \left( |\nabla E|^2 - {1\over 12}
    |\nabla R|^2 \right) + 6tr E^3 + R|E|^2 -
      6W_{ijk\ell} E_{ik} E_{j \ell}\right\} dvol \,,
\tag 2.6
$$
where $tr E^3 = E_{ij} E_{ik} E_{jk}$, and
$$
\int\limits_{M^4} |\nabla W |^2 dvol =
  \int \left\{ 72 \det W^+ + 72 \det
W^- -
  {1\over 2} R |W|^2 + 2 W
_{ij\kappa\ell}E_{i\kappa} E_{j\ell} \right\}dvol, \tag 2.7
$$
where $W^+$ and $W^-$ are respectively the self-dual and the
anti-self dual part of the Weyl tensor.
\endproclaim

We remark that the second identity is a consequence of Stokes'
theorem, the Bianchi identities, and the definition of the Bach
tensor in (2.4). In this paper we do not need the precise form of
the identity (2.7) but only the fact that the terms in $\det W^+$
and $\det W^-$ are of the form $W*W*W$ --a contraction of three Weyl
tensor terms.

In the recent work of Tian-Viaclovsky--which is obtained from an
iteration process by applying equation (2.5) and the the Sobolev
inequality (2.2)- they have derived the following "$\epsilon$
-regularity" theorem.

\proclaim{Theorem 2.4}([TV-1], Theorem 3.1) \quad Suppose that
$(M^4, g)$ is a Bach flat 4-manifold with the Yamabe constant $Y(
M^4, [g]) >0$ and let $g$ denotes its Yamabe metric. Then there
exist some positive numbers $\tau_k$ and $C_k$ depending on $Y (M^4,
[g]) $ such that, for each geodesic ball $B_{2r}(p)$ centered at
$p\in M$, if
$$
\int_{B_{2r}(p)} |Rm|^2 dv \leq \tau_k, \tag 2.8
$$
then
$$
\sup_{B_r(p)} |\nabla^k Rm| \leq  \frac {C_k}{r^{2+k}}
(\int_{B_{2r}(p)} |Rm|^2 dv)^\frac 12. \tag 2.9
$$
\endproclaim

Built on the estimate in Theorem 2.4, the main result in [TV-1]
--which plays a crucial role in this paper -- is the following
volume estimates (2.11) of geodesic balls:

\proclaim{Theorem 2.5} Let $(X, g)$ be a complete, non-compact, four
manifold with base point $p$, and let $r(x) = d (p, x),$ for $ x \in
X.$ Assume that there exists a constant $v_0 >0$ so that
$$
\text{vol} (B_r(q)) \geq v_0 r^4 \tag 2.3
$$
holds for all $ q \in X$,  assume furthermore that as $r \rightarrow
\infty$,
$$
sup_{S(r)} |Rm_g| = o ( r^{-2} ), \tag 2.10
$$
where $S(r) = \partial B_r (p).$  Assume further that the first
Betti number $b_1 (X)  < \infty,$ then $(X,g)$ is an ALE space, and
there exists a constant $v_1$ (depending on $(X, g)$) so that
$$
\text{vol} (B_r(p)) \leq v_1 r^4.  \tag 2.11
$$
\endproclaim

We remark that the constant $v_1$ in (2.11) obtained above may
depend on the given manifold $(X,g)$.

We return to the class of metrics $g$ with positive Yamabe constant
and which satisfy the inequality (1.10) in the statement of Theorem
B. First, we make some simple observations:

\proclaim{Lemma 2.6} Suppose that $(M^4, g)$ is a closed, compact
4-manifold with positive Yamabe constant and which satisfies the
inequality (1.10). Suppose further that $g$ is a Yamabe metric with
$vol (g) = \frac {8 \pi^2}{3}$, then
$$
\frac {3}{\pi} \sqrt a_0 \leq R [g] \leq 12, \tag 2.12
$$
and
$$
\int_M |E|^2 dv_g \leq (32\pi^2 - 2 a_0). \tag 2.13
$$
\endproclaim

\vskip .1in We remark that if the metric $g$ satisfies the
inequality (1.9) in the statement of the gap Theorem A, then $ a_0 =
(1- \epsilon) 16 \pi^2$ and we have
$$
12 \sqrt {(1- \epsilon)} \leq R [g] \leq 12. \tag 2.12'
$$
and
$$
\int_M |E|^2 dv_g \leq 32\pi^2 \epsilon. \tag 2.13'
$$

A deeper result for this class of metrics is following result of
Gursky:

\proclaim{Theorem 2.7}([G]) Suppose that $(M^4, g)$ is a closed
4-manifold with positive Yamabe constant and $\int_M \sigma_2(A_g)
dv_g > 0.$  Then $ b_1(M^4) = 0 $.
\endproclaim

We remark that under the same assumptions as in Theorem 2.7, in
[CGY-1], we obtained the stronger result that the fundamental group
of $M^4$ is finite.

In summary, we have:

\proclaim{Corollary 2.8} Suppose that $\Cal A$ is a collection of
Bach flat Riemannian manifolds $(M^4, g)$ with positive Yamabe
constants, satisfying (1.8) and either (1.9) or (1.10), then there
is a lower bound of the Yamabe constants and an upper bound of the
$L^2$ norm of the curvature tensor $Rm_g$ for the whole collection;
also $b_1 (M^4) = 0$. Thus one may apply results in Theorem 2.4 and
and Theorem 2.5 to metrics in the collection $\Cal A$.
\endproclaim

\vskip 0.1in \noindent{\bf 3. The Neck Theorem}\vskip 0.1in

The neck theorem which we are going to establish in this section is
modeled after the result of Anderson-Cheeger (Theorem 1.18 in [AC]),
in which the theorem was established under the assumption that there
is a point-wise Ricci curvature bound on the manifold.

%% should we add some sentences here?

Suppose $(M^4, g)$ is a Riemannian manifold. For a point $p\in M$,
we denote $B_r(p)$ the geodesic ball with radius $r$ centered at
$p$, $S_r(p)$ the geodesic sphere of radius $r$ centered at $p$. We
consider the geodesic annulus centered at $p$ as:
$$
\bar A_{r_1, r_2} (p) = \{ q\in M: r_1 \leq \text{dist}(q, p) \leq
r_2\}. \tag 3.1
$$
In general, $\bar A_{r_1, r_2}(p)$ may have more than one connected
components. We will denote by
$$
A_{r_1,r_2}(p) \subset \bar A_{r_1, r_2}(p)
$$
any component of $\bar A_{r_1,r_2}(p)$ which meets $S_{r_2}(p)$:
$$
A_{r_1,r_2}(p)\bigcap S_{r_2}(p) \neq \emptyset. \tag 3.2
$$
Since the distance function is Lipschitz, we may consider the 3
dimensional Hausdorff measure for the geodesic sphere $S_r(p)$ and
denote it by $\Cal H^3(S_r(p)).$

\proclaim{Theorem 3.1} Suppose that $(M, g)$ is a Bach flat
Riemannian 4-manifold with a Yamabe metric with a positive Yamabe
constant, and suppose the first Betti number $b_1 (M) = 0$. Assume
that $p \in M$, and $\alpha \in (0,1)$, $\epsilon >0$, $v_2 >0$  and
$a < \text{dist}(p,
\partial M)$ are given constants.
Then there exist positive numbers $\delta_0, c_2, n$ depending on
$\epsilon,\, \alpha,\, C_s,\, v_2 \, \, and \,\,  a \, $ such that
the following statements hold. Let $A_{r_1, r_2}(p)$ be a connected
component of the geodesic annulus which satisfies the condition
(3.2), with
$$
r_2 \leq c_2 a, \quad\quad r_1 \leq \delta_0 r_2, \tag 3.3
$$
$$
\Cal H^3 (S_r(p)) \leq v_2 r^3, \quad\forall r\in [r_1, 100r_1],
\tag 3.4
$$
and
$$
\int_{A_{r_1, r_2}(p)} |Rm|^2 dv \leq \delta_0. \tag 3.5
$$
Then $A_{r_1, r_2}(p)$ is the only component which satisfies (3.2).
In addition for this unique component
$$
A_{(\delta_0^{-\frac 14}-\epsilon)r_1, (\delta_0^\frac 14 +\epsilon)
r_2}(p),
$$
which intersects with $S_{(\delta_0^\frac 14 +\epsilon) r_2}(p)$,
there exist some finite group $\Gamma \subset O(4)$, acting freely
on $S^3$, with $|\Gamma| \leq n$, and a quasi isometry $\Psi$, so
that
$$
A_{(\delta_0^{-\frac 14}+\epsilon)r_1, (\delta_0^\frac 14 -\epsilon)
r_2}(p) \subset \Psi(C_{\delta_0^{-\frac 14}r_1, \delta_0^\frac 14
r_2} (S^3/\Gamma)) \subset A_{(\delta_0^{-\frac 14}-\epsilon)r_1,
(\delta_0^\frac 14 +\epsilon) r_2}(p)
$$
such that for all $C_{\frac 12 r, r}(S^3/\Gamma) \subset
C_{\delta_0^{-\frac 14}r_1, \delta_0^\frac 14 r_2} (S^3/\Gamma)$, in
a local coordinate, one has
$$
|(\Psi^*(r^{-2}g))_{ij} - \delta_{ij}|_{C^{1, \alpha}} \leq
\epsilon. \tag 3.6
$$
\endproclaim

We now begin the proof of Theorem 3.1; the first part of our
arguments are in large part modifications of the arguments used in
the proof of Theorem 4.1 in [TV-1] (see Theorem 2.5 for the
statement of the theorem); the second part some modification of the
proof of Theorem 1.18 in [AC]. We will be brief in presenting our
proof here. We begin with a lemma to prove the uniqueness of the
connected component satisfying condition (3.2).

\proclaim{Lemma 3.2} Under the same assumptions of the theorem,
there is only one connected component in the annulus $\bar A_{r_1,
r_2}(p)$ which intersects with $S_{r_2}(p)$.
\endproclaim

\demo{Proof} Assume otherwise, that is, assume there is another
component $B_{r_1, r_2}(p)$ in $\bar A_{r_1, r_2}(p)$ which also
intersects with $S_{r_2}(p)$. It then follows from the assumption
$b_1 (M) =0 $ that these two components can not be joined by a curve
from their intersections with the components of $S_{r_2}(p)$. Thus
the region $A_{r_1, 2r_1}(p) \subset A_{r_1, r_2}(p)$ separates $M -
A_{r_1, 2r_1}(p)$ into two disjoint components, so that the volume
of each components at least $C r_2^4$ for some constant C. We now
observe that under the assumption that the Yamabe constant $Y= Y(M,
[g])$ be positive, the Sobolev inequality (2.1) holds by Lemma 2.1.
If we define a function $u$ to be equal to $1$ on the component with
a smaller volume and $0$ on the other component and set it to be a
smooth function which varies from $0$ to $1$ with its first
derivative bounded by $ \frac {2}{r_1}$ on $A_{r_1, 2r_2}(p)$, then
apply the Sobolev inequality (2.1) and the assumption (3.4) we have
$$
\aligned Y (C r_2^4)^\frac 12  & \leq Y (\int_M u^4dv)^\frac 12 \leq
6 \int_M |\nabla u|^2 dv + \frac {Y} {\sqrt{\text{vol}(M)}} \int_M
u^2dv \\
& \leq \frac {24}{r_1^2}v_2(16r_1^4 - r_1^4) + \frac {Y}
{\sqrt{\text{vol}(M)}} v_2 r_1^4 + \frac {Y} {\sqrt 2}(\int_M
u^4dv)^\frac 12.
\endaligned
\tag 3.7
$$
This implies
$$
r_2^4 \leq  C' r_1^4, \tag 3.8
$$
for some suitable constant $C'$ and contradicts with the assumption
(3.3) of the Lemma.
\enddemo

\demo{Proof of Theorem 3.1} \vskip .1in

\noindent $\underline {Part \, \, I  \,:}$ \quad   First we apply
Theorem 2.4 in section 2 to obtain
$$
\sup_{S_r(p)\bigcap A_{r_1, r_2}(p)} |\nabla^k Rm| \leq \frac
{C_k\epsilon (r)}{r^{2+k}} \delta_0^\frac 12, \tag 3.9
$$
for each $r \in [2r_1, \frac 12 r_2]$, provided that $\delta_0 \leq
\tau_k$ for $k = 0, 1$, where $\epsilon (r)\leq 1$ and goes to $0$
as $r/r_1\to\infty$.

We then observe that by arguments similar to that in the proof of
Lemma 3.2 above, there is a unique connected component which
intersects with $S_{\frac 12 r_2}(p)$, which we denote by $A_{2r_1,
\frac 12 r_2}(p)$. This component $A_{2r_1, \frac 12r_2}(p)$ has a
tree structure in the sense that while the intersections of $S_s(p)$
with the connected components $A_{r,s}(p)$ may have more than one
connected components; the intersections of $S_r(p)$ with the
connected components $A_{r, s}(p)$ are always connected for all
$2r_1 < r< s< \frac 12 r_2 $. Following the argument in section 4 in
[TV-1], we choose a constant $s$ with $2\leq s \leq 10$, and denote
by $N$ the integer determined by the condition $r_2\in [2s^Nr_1,
2s^{N+1}r_1)$. Let $\{A_{2s^jr_1, 2s^{j+1}r_1}(p)\}_{j=0}^N $ denote
the set of annuli such that
$$
S_{2s^{j+1}r_1}(p)\bigcap A_{2s^{j+1}r_1, 2s^{j+2}r_1}(p)\subset
S_{2s^{j+1}r_1}(p)\bigcap A_{2s^{j}r_1, 2s^{j+1}r_1}(p).
$$
and call $$ D(s) = \bigcup_{j=0}^{N}A_{2s^{j}r_1, 2s^{j+1}r_1}(p)
\tag 3.10
$$
a direction in the tree $A_{2r_1, \frac 12 r_2}(p)$.  We claim that
there exists some constants $\delta_1>0$ and $ v_3>0$ such that
$$
\Cal H^3 (S_r(p)\bigcap D(s)) \leq v_3 r^3, \quad for \quad  all
\quad  r \in [2r_1, \frac 14 r_2], \tag 3.11
$$
provided that $\delta_0 \leq \delta_1$. Note that as a consequence
of (3.11), we  have
$$
\aligned \Cal H^3 (S_r(p) \bigcap A_{2s^jr_1, 2 s^{j+1}r_1}(p)) &
\leq  C_4 r^3 \\ \text{vol}(A_{2s^{j}r_1, 2s^{j+1}r_1}(p)) & \leq
C_4 s^{4j}r_1^4 \endaligned \tag 3.12
$$
for some constant $C_4$.

The proof to establish the claim (3.11) in [TV-1] is rather
complicated. One of the key step there is to show that given a
direction, there exists some maximal subsequence of annuli
$A_{2s^{j}r_1, 2 s^{j+1}r_1}(p)$, such that
$$
\Cal H^3 (S_{2s^{j+1}r_1}(p)\bigcap A_{2s^{j+1}r_1, 2s^{j+2}r_1}(p))
\geq (1-\eta_j) \Cal H^3(S_{2s^{j}r_1}(p)\bigcap A_{2s^{j}r_1,
2s^{j+1}r_1}(p)) \tag 3.13
$$
for some sequence $\eta_j \to 0$. One then establishes the estimates
in (3.12) using a proof by contradiction for this subsequence of
annuli. We observe that the same strategy of proof works in our
case. That is, assuming that (3.11) does not hold, then there is a
sequence of of $(M_i, g_i)$ and directions in the annuli
$$
D^i(s) \subset A^i_{2r_1^i, \frac 12 r^i_2}(p_i) \subset M_i,
$$
satisfying (3.3), (3.4) and (3.5) for some sequence $r_1^i, r_2^i
\to 0$, $\delta^i \to 0$, and a maximal subsequence of annuli
$A^i_{2s^jr_1^i, 2 s^{j+1}r_1^i}(p_i)$ and numbers $\eta_{i, j} \to
0$ such that
$$
\aligned \Cal H^3 & (S^i_{2s^{j+1}r_1^i}(p_i)\bigcap
A^i_{2s^{j+1}r^i_1, 2s^{j+2}r^i_1}(p_i)) \\ & \geq (1-\eta_{i,j})
\Cal H^3(S^i_{2s^{j}r^i_1}(p_i)\bigcap A^i_{2s^{j}r^i_1,
2s^{j+1}r^i_1}(p_i))
\endaligned
$$
due to the assumption (3.4) in the statement in Theorem 3.1. So that
$$
\text{vol}(A^i_{2s^{j}r^i_1, 2s^{j+1}r^i_1}(p_i)) \to \infty
$$
when both $i$ and $j$ tend to infinity. We then follow the same line
of argument as in [TV-1] to get a contradiction and establish the
claims (3.11) and (3.12).

We may now repeat the argument similar to that used in the proof of
Lemma 3.2 above to show that the direction $D(s)$ established does
not further split into more branches at the end sufficient far away
from the initial sphere $S_{r}(p)$ of the annulus. More precisely:

\roster
\item"{(3.14)}"  There is a constant $K$,
which only depends on the Sobolev constant in (2.2),  such that
there is only one component in the geodesic annulus $\bar
A_{r,t}(p)$ which intersects with $S_{t}(p)$, when $2r_1 \leq r< Kr
\leq t \leq \frac 1{2s} r_2$. \vskip .1in
\item"{(3.15)}"  Moreover the estimates in (3.12) hold for all
geodesic annuli in between $2r_1$ and $\frac 1{2s} r_2$.
\endroster

\vskip 0.1in\noindent $\underline {Part \, \,  II \, :}$ \quad We
will now modify the argument in the proof of Theorem 1.18 in [AC] to
finish the rest of the proof of Theorem 3.1.

Again, we argue by contradiction. Suppose the statement of the
theorem is not true, then for some $\epsilon_0$, there exist a
sequence of $(M_i, g_i)$ and annuli
$$
A^i_{r_1^i, r^i_2}(p_i) \subset M_i
$$
satisfying (3.3), (3.4) and (3.5) for some sequence $r_1^i, r_2^i
\to 0$, and $\delta^i \to 0$, such that some sub-annuli $A^i_{\frac
12 r_i, r_i}(p_i)$ with $r_i\in [(\delta^i)^{-\frac 14}r_1^i,
(\delta^i)^{\frac 14}r_2^i]$ with rescaled metric $\tilde g_i =
r_i^{-2}g_i$, are not $\epsilon_0$-close in $C^{1, \alpha}$ topology
to annuli in any cone $C(S^3/\Gamma)$. We note that by our choice of
$r_i$, we have $\frac {r_2^i}{r_i} \to \infty,$ thus there is no
curvature concentration on the annuli $A^i ( \frac 1l r_i, l r_i)$
for any positive integer $l$.

Consider the sequence $(M_i, r_i^{-2}g_i, q_i)$ of manifolds with
base points $q_i \in S_{r_i} (p_i). $  For each integer $l$, denote
$A_{\frac 1l r_i, l r_i}(p_i)$  a connected component in the
geodesic annuli \quad $\bar A_{\frac 1l r_i, l r_i}(p_i)$ which
reaches to $S_{l r_i}(p_i)$; by (3.14) in Part I, such a component
$A_{\frac 1l r_i, l r_i}(p_i)$ is unique.

Now fix $l$ and consider the spaces $(A_{\frac {l} {r_i}, l
r_i}(p_i), r_i^{-2}g_i, q_i)$. We first assert that it follows from
our assumption of the uniform positive lower bound of the Yamabe
constants on the manifolds $M_i$, that we may apply Lemma 2.2 to
have some lower bound which is uniform in scale for the volume of
the geodesic balls (see (2.3)). By applying condition (3.15) in Part
I above, we also have a bound for the {\it intrinsic} diameter of
the boundary $S_{\frac 1l r_i} (p_i)$ of $ A_{\frac 1l r_i, l
r_i}(p_i)$ which is again uniform in scale. Thus (see [An-1]) the
spaces converge in the Gromov-Hausdorff topology as $i \to \infty $
to a space which we denote by $(B_\infty^l, g_\infty^l)$; note that
we have $(B_\infty^l, g_\infty^l) \subset (B_\infty^{l+1},
g_\infty^{l+1}).$

Denote
$$
(B_\infty, g_\infty) = \bigcup_{l = k}^\infty (B^l_\infty,
g^l_\infty).
$$
Again, (3.15) in Part I above implies that $p_{\infty}$ is an
isolated singularity for the space $(B_\infty, g_\infty).$  Thus it
follows from the argument of Theorem 1.18 in [AC] that the limiting
space $B_{\infty}$ is a Euclidean cone. Finally we observe that by
(3.14), the geodesic annulus $A_{\frac {r_i}{\sqrt K},  \sqrt K r_i}
(p_i)$ is contained in $A_{\frac 1l r_i, l r_i} (p_i)$ for $l$
sufficient large thus it tends to part of the Euclidean cone. We can
then piece-wisely connect such annuli and conclude the statement in
the theorem by the same argument as in ([AC], page 241).

We have thus finished the proof of Theorem 3.1.

\vskip 0.1in \noindent{\bf 4. Bubble tree construction}\vskip 0.1in

In this section we consider a sequence of metrics $(M_i, g_i)$
satisfying the conditions stated in the next paragraph. We study the
metrics near points of curvature concentration. When there is
concentration of curvature, we will blow up the metrics near certain
prescribed points in order to extract convergent subsequence of
rescaled metrics. The main tool used is the Neck Theorem of the
previous section. We follow the procedure used in [BC], [Q] [St] to
extract bubbles, proceeding from the smaller scales to larger scales
and ending at the the largest bubble. This procedure is different
from that of Anderson and Cheeger ([AC]), whch starts from the
larger scales to the smaller scales. We find it necessary to proceed
in this manner because our Neck theorem is weaker than that of [AC],
since we do not have a priori diameter control for the bubble body
while in the case considered by [AC] such control is implied by the
Ricci curvature bounds.

We will assume in this section that $(M_i, g_i)$ satisfy the
following conditions:

\vskip .1in

\noindent (4.1) $(M_i.g_i)$ are Bach flat closed 4-manifolds,

\noindent (4.2) There is a positive lower bound for the Yamabe
invariants:

$$
Y(M_i, [g_i]) = \inf_{g\in [g_i]} \frac {\int_{M_i} R_g
dv_g}{(\int_{M_i}dv_g)^\frac 12} \geq Y_0 >0,
$$
for some fixed number $Y_0$,

\noindent (4.3) $b_1(M_i)=0$,

\noindent (4.4) There is a common  bound for the curvature tensor:
$$
\int_{M_i} |Rm^i|^2 dv^i \leq \Lambda,
$$
for some fixed number $\Lambda$.

In the following we will ignore the difference between a sequence
and its subsequences for simplicity, since it will not cause any
problem for the arguments.

Let us start the construction of bubble tree. Set $\delta = \min\{
\tau_0, \tau_1, \delta_0\}$ where $\tau_k$ are from the
$\epsilon$-estimates for curvature (cf. Theorem 2.5) and $\delta_0$
is from the neck Theorem 3.1. This will be an iterative
construction.
% so we suppose that $X_i \subset M_i$
%is a region containing  a geodesic ball of a fixed radius
%$r_0$ and its boundary does not carry curvature, i.e.
%$$
%\int_{T_{\eta_0}(\partial X_i)}
%|Rm^i|^2 dv^i \leq \frac \delta 2,
%$$
%where
%$$
%T_{\eta_0}(\partial X_i) = \{p\in M_i: \text{dist}(p,
%\partial X_i) < \eta_0 \},
%$$
%for some fixed positive number $4\eta_0 < r_0$.
For $ p \in M_i $, denote $B_r^i (p) $ a geodesic ball of radius $r$
centered at $p$ in $(M_i, g_i)$,
$$
s_i^1(p) = r \ \text{such that} \ \int_{B^i_r(p)} |Rm^i|^2dv^i =
\frac \delta 2. \tag 4.5
$$
We then choose
$$
p_i^1 = p  \ \text{such that} \ s^1_i(p) = \inf_{p \in M_i} s_i^1
(p). \tag 4.6
$$
We may assume that $\lambda_i^1 = s_i^1(p_i^1)\to 0$, for otherwise
there would be no curvature concentration in $M_i$. We then conclude
that $(M_i, (\lambda_i^1)^{-2}g_i, p_i^1)$ converges to
$(M_\infty^1, g_\infty^1, p_\infty^1)$, which is a Bach flat, scalar
flat, complete 4-manifold satisfying the Sobolov inequality (2.2),
and $L^2$ bound for the Riemann curvature (4.4), and in addition
having only one end.

\proclaim{Definition 4.1} We call a Bach flat, scalar flat, complete
smooth 4-manifold with the Sobolev inequality (2.2),
$L^2$-Riemannian curvature (4.4), and a single ALE end a {\it leaf
bubble}, while we will call such a space with finitely many isolated
irreducible orbifold points an {\it intermediate bubble}.
\endproclaim

It follows that on such a bubble, the volume growth tends to that of
the standard Euclidean cone $C(S^3/\Gamma_1)$ for some $\Gamma_1$
near the leaf bubble on $(M_\infty^1, g_\infty^1)$. Therefore, by
Theorem 2.5 for some $K^1$ (which is allowed to depend on
$(M_\infty^1, g_\infty^1)$ at this point in our approach) such that
$$
\int_{M_\infty^1 \setminus B_{\frac {K^1}2}(p_\infty^1)} |Rm|^2 d
v_{\infty}^1 \leq \frac \delta 4 \tag 4.7
$$
and
$$
\Cal H^3(S_r(p_\infty^1)) \leq 2 v_4 r^3, \tag 4.8
$$
for all $S_r(p_\infty^1) \subset M_\infty^1$ and $r\in [K^1,
100K^1]$, where $v_4$ may be chosen to be the one for Euclidean
space $(R^4, |dx|^2)$ and $4v_4 = v_2$ in Theorem 3.1.

If there is further curvature concentration, we proceed to extract
the next bubble. We define, for $p\in M_i \setminus
B^i_{K^1\lambda_i^1} (p_i^1)$,
$$
s_i^2(p) = r \ \text{such that} \ \int_{B_r^i (p) \setminus
B^i_{K^1\lambda_i^1} (p_i^1)} |Rm^i|^2dv^i = \frac \delta 2. \tag
4.9
$$
We then choose
$$
p_i^2 = p  \ \text{such that} \ s^2_i(p) = \inf_{M_i
 \setminus B^i_{K^1\lambda_i^1} (p_i^1)} s_i^2(p). \tag 4.10
$$
Again we may assume that $\lambda_i^2 = s_i^2(p_i^2)\to 0$, for if
otherwise there would be no more curvature concentration. It follows
from (4.7), (cf. (58) in [Q]) that:

\proclaim{Lemma 4.2}
$$
\frac {\lambda_i^2}{\lambda_i^1} + \frac {\text{dist}(p_i^1,
p_i^2)}{\lambda_i^1} \to \infty. \tag 4.11
$$
\endproclaim

\demo{Proof} Assume to the contrary, there exists a number $M>0$
such that
$$
1 \leq \frac {\lambda_i^2}{\lambda_i^1} \leq M \ \text{and} \ \frac
{\text{dist}(p_i^1, p_i^2)}{\lambda_i^1} \leq M.
$$
Then
$$
p_i^2 \in B^i_{M\lambda_i^1}(p_i^1).
$$
Hence
$$
B^i_{\lambda_i^2}(p_i^2)\setminus B^i_{K^1\lambda_i^1}(p_i^1)
\subset B^i_{3M\lambda_i^1}(p_i^1) \setminus
B^i_{K^1\lambda_i^1}(p_i^1)
$$
and this contradicts (4.9). This finishes the proof of the Lemma.
\enddemo

There are two possible cases which we will discuss separately:
$$
\aligned & \text{Case 1.} \quad  \frac {\text{dist}(p_i^1,
p_i^2)}{\lambda_i^2} \to \infty \, \,\,  as \,\, i \to \, \infty,  \\
& \text{Case 2.} \quad \frac {\text{dist}(p_i^1,
p_i^2)}{\lambda_i^2} \leq M^1, \, \, for \,\, some \,\,  constant
\,\, M^1 \,\, as \,\,  i \,\,  \to \, \infty.
  \endaligned \tag 4.12
$$
In case 1, applying Lemma (4.2),  we certainly also have
$$
\frac {\text{dist}(p_i^1, p_i^2)}{\lambda_i^1} \to \infty.
$$
Therefore, in the limit of the convergent sequence of metrics $(M_i,
(\lambda_i^2)^{-2}g_i, p_i^2)$ one does not see the concentration
which produces the bubble $(M_\infty^1, g_\infty^1)$, and converges
to another end bubble $(M_\infty^2, g_\infty^2)$. Similarly in the
limit of the convergent sequence of metrics $(M_i,
(\lambda_i^1)^{-2}g_i, p_i^1)$ one does not see the concentration
which produces $(M_\infty^2, g_\infty^2)$. It follows from (4.5)
that there are at most $2\Lambda/\delta$ number of such leaf
bubbles.

\proclaim{Definition 4.3} We say two bubbles $(M_\infty^{j_1},
g_\infty^{j_1})$ and $(M_\infty^{j_2}, g_\infty^{j_2})$ associated
with $(p_i^{j_1}, \lambda_i^{j_1})$ and $(p_i^{j_2},
\lambda_i^{j_2})$ are {\it separable} if
$$
\frac{\text{dist}(p_i^{j_1}, p_i^{j_2})}{\lambda_i^{j_1}} \to \infty
\ \text{and} \ \frac{\text{dist}(p_i^{j_1},
p_i^{j_2})}{\lambda_i^{j_2}} \to \infty. \tag 4.13
$$
\endproclaim
In Case 2, one starts to introduce intermediate bubbles which will
contain in its body a number of previously constructed bubbles. It
would be helpful to observe that, in contrast to [BC] [Q] [St], one
needs the neck Theorem 3.1 to produce the intermediate bubbles. If
the Ricci curvature is bounded as in the case considered in [AC],
one would have no problem to take limit in Gromov-Hausdorff topology
to be a complete 4-manifold with finitely many possible orbifold
points. Instead we will apply our neck Theorem 3.1 to prove that the
limit space has only isolated point singularities, which are then
orbifold points according to [TV-2].

\proclaim{Lemma 4.4} Suppose that there are several separable
bubbles $\{(M_\infty^j, g_\infty^j)\}_{j\in J}$ associated with
$\{(p_i^j, \lambda_i^j)\}_{j\in J}$,and  there is further
concentration detected as $(p_i^k, \lambda_i^k)$ after $\{(p_i^j,
\lambda_i^j)\}_{j\in J}$, such that
$$
\frac {\text{dist}(p_i^j, p_i^k)}{\lambda_i^k} \leq M^j, \tag 4.14
$$
for some constant $M^j$ as $i$ tends to $\infty$; and hence,
$$
\frac {\lambda_i^k}{\lambda_i^j} \to \infty
$$
for each $j\in J$. Suppose, in addition, that $\{(p_i^j,
\lambda_i^j)\}_{j\in J}$ is a maximal collection of such bubbles,
then $(M_i, (\lambda_i^k)^{-2}g_i, p_i^k)$ converges in
Gromov-Hausdorff topology to an intermediate bubble $(M_\infty^k,
g_\infty^k)$. We will call such a bubble a parent or grandparent of
the given collection of bubbles $\{(M_{\infty}^j,
g_{\infty}^j)\}_{j\in J}$.\endproclaim

\demo{Proof} It follows from Theorem 2.5, that there are constant
$K^j$ for each bubble such that
$$
\int_{M_\infty^j \setminus B_{\frac {K^j}2}(p_\infty^j)} |Rm|^2 d
v^j< \frac \delta 4 \tag 4.15
$$
and
$$
\Cal H^3(S_r(p_\infty^j)) \leq 2 v_4 r^3, \tag 4.16
$$
for all $S_r(p_\infty^j) \subset M_\infty^j$ and $r\in [K^j,
100K^j]$. The claim is clear if $J$ consists of only a single
element. If $J$ contains more than one element, we proceed in two
steps. \vskip 0.1in \noindent {\it Step 1:} \quad We first prove the
lemma in the case that there exists some constant $\eta_k > 0$ so
that
$$
\frac {\text{dist}(p_i^j, p_i^{j'})}{\lambda_i^k} \geq \eta_k > 0,
\,\, for \,\,  all\,\, j, \, j' \,  \in \,\,  J. \tag 4.17
$$

For any given number $L>>\eta_k >>1/L$ and $i>>1$ with $L \geq M^j $
for each $ j \in J$,  we have
$$
B^i_{\lambda_i^k L} (p_i^k) \supset B^i_{\lambda_i^k /L}(p_i^j)
\supset B^i_{\lambda_i^j K^j}(p_i^j) \tag 4.18
$$
for each $j\in J$. It follows from (4.15) and (4.16) that we may
take limit in rescaled sequence:
$$
(B^i_{\lambda_i^k L} (p_i^k)\setminus (\bigcup_j B^i_{\lambda_i^k
/L}(p_i^j)), (\lambda_i^k)^{-2}g_i, p_i^k) \to (M_\infty^k(L),
g_\infty^k(L)) \tag 4.19
$$
for each given $L$, and
$$
(M_\infty^k(L), g_\infty^k(L)) \subset (M_\infty^k(L+1),
g_\infty^k(L+1)) \subset \bigcup_{l=1}^\infty (M_\infty^k(L+l),
g_\infty^k(L+l)). \tag 4.20
$$
We will now apply the neck theorem to show that
$$
(M_\infty^k, g_\infty^k) = \bigcup_{l=1}^\infty (M_\infty^k(L+l),
g_\infty^k(L+l))
$$
is an intermediate bubble. It follows from (4.17) and the choice of
$\lambda_i^k$ that (4.18) holds and
%$$
%B^i_{\lambda_i^j K^j}(p_i^j) \subset B^i_{\lambda_i^k \frac
%1L}(p_i^j) \subset B^i_{\lambda_i^k}(p_i^k) \tag 4.21
%$$
$$
\int_{B^i_{\lambda_i^k \frac 1L}(p_i^j)\setminus B^i_{\lambda_i^j
K^j}(p_i^j)}|Rm^i|^2 dv^i \leq \frac \delta 2 < \delta_0. \tag 4.21
$$
The Neck theorem then shows that the diameter of $S^i_{\lambda_i^k
\frac 1L}(p_i^j)$ in the rescaled space $(B^i_{\lambda_i^k L}
(p_i^k)\setminus (\bigcup_j B^i_{\lambda_i^k /L}(p_i^j)),
(\lambda_i^k)^{-2}g_i, p_i^k)$ goes to zero.

\vskip 0.1in \noindent {\it Step 2:} \quad  We now deal with the
case that there is a subset $J'\subset J$ such that
$$
\frac {\text{dist}(p_i^j, p_i^{j'})}{\lambda_i^k} \to 0 \tag 4.22
$$
for all $j, j'\in J'$ as $i \to \infty$. Our strategy will be to
combine some elements in $J'$ to create some intermediate bubbles.
We remark that this situation does not arise in the case considered
in [AC] due to the gap theorem for Ricci flat complete orbifolds
(cf. [Ba]). We call those intermediate bubbles thus created the
exotic bubbles, since they may carry an arbitrarily small amount of
$L^2$-Riemannian curvature.  To start the creation process, let
$$
\mu_i^1 = \min\{ \text{dist}(p_i^j, p_i^{j'}): j, j'\in J'\} =
\text{dist}(p_i^{j_1}, p_i^{j_2}), \tag 4.23
$$
$$
J_1 = \{j\in J': sup_i\{\frac {\text{dist}(p_i^{j_1},
p_i^j)}{\mu_i^1}\} < \infty\}, \tag 4.24
$$
note that $J_1$ consists of at least two elements $j_1$ and $j_2$;
and for some large number $N_1$
$$
\frac {\text{dist}(p_i^{j_1}, p_i^j)}{\mu_i^1} \leq \frac 12 N_1.
\tag 4.25
$$
Then we consider the sequence $(M_i, (\mu_i^1)^{-2}g_i, p_i^{j_1})$.
In view of Theorem 2.5 and the volume bound we may take limit in the
rescaled sequence:
$$
(B^i_{(N_1+l)\mu_i^1}(p_i^{j_1})\setminus (\bigcup_{j\in J_1}
B^i_{\frac 1{N_1+l}\mu_i^1}(p_i^j)), (\mu_i^1)^{-2}g_i, p_i^{j_1})
\,\, \to (N_{\infty}^1(l),h_{\infty}^1(l)).$$ Then
$$
(N_\infty^1, h_\infty^1) = \bigcup_{l=1}^\infty (N_\infty^1(l),
h_\infty^1(l))
$$
gives an exotic intermediate bubble by applying the neck Theorem 3.1
to the annuli $B^i_{\frac 1{N_1}\mu_i^1}(p_i^j)) \setminus
B^i_{K^j\lambda_i^j}(p_i^j)$ for all $j\in J_1$. This exotic bubble
$(N_\infty^1, h_\infty^1)$ has at most $|J_1|$ number of orbifold
singularities. We use this bubble to replace all bubbles
$(M_\infty^j, g_\infty^j)$ for $j\in J_1$. In other words, we
combine all bubbles $(M_\infty^j, g_\infty^j)$ for $j\in J_1$ into
one single bubble $(N_\infty^1, h_\infty^1)$. Since
$$
\frac {\text{dist}(p_i^{j_1}, p_i^j)}{\mu_i^1} \to \infty
$$
for all $j\in J'\setminus J_1$ by the definition of $J_1$ and
$$
\frac {\text{dist}(p_i^{j_1}, p_i^j)}{\lambda_i^j} \to \infty
$$
as assumed, the collection
$$
\{(p_i^{j_1}, \mu_i^1)\}\bigcup \{(p_i^j, \lambda_i^j): j\in
J'\setminus J_1\}
$$
constitute a smaller family of separable bubbles. We may then repeat
the above combination process for a finite number of steps if
necessary, to combine all bubbles $\{(M_\infty^j,
g_\infty^j)\}_{j\in J'}$ into one single bubble $(N_\infty,
h_\infty)$ associated with some $(p_i^l, \mu_i^l)$ with $l \in J'$
and
$$
\frac {\lambda_i^k}{\mu_i^l} \to \infty \ \text{and} \ \frac
{\text{dist}(p_i^l, p_i^k)}{\lambda_i^k} \leq M^l. \tag 4.26
$$
Let us call the new collection
$$
\{(p_i^{j\prime}, \mu_i^{j\prime})\}\bigcup \{(p_i^j, \lambda_i^j):
j\in J \setminus J^{\prime}.\} \tag 4.27
$$
Now we are back to the situation in which the new collection of
bubbles (4.27) satisfy hypothesis (4.14) of Lemma 4.4. Thus we may
repeat the process of step 1 at the beginning of the proof of the
lemma to the new collection and note that after finite many such
repetitions, the final collection of bubbles will satisfy (4.17). We
have thus finished the proof of Lemma 4.4.
\enddemo

\proclaim{Definition 4.5} A {\it bubble tree} $T$ is defined to be a
tree whose vertexes are bubbles and whose edges are necks from
Theorem 3.1. At each vertex $(M_\infty^j, g_\infty^j)$, its ALE end
is connected, via a neck, to its parent toward the {\it root bubble}
of $T$, while at finitely many isolated possible orbifold points of
$(M_\infty^j, g_\infty^j)$, it is connected, via necks, to its
children toward {\it leaf bubbles} of $T$. We say two bubble trees
$T_1$ and $T_2$ are {\it separable} if their root bubbles are
separable.
\endproclaim

To finish the bubble tree construction, we describe the following
inductive procedure. Suppose that there are $m$ bubbles in $n (\leq
m)$ separable bubble trees $\{T^j\}_{j=1}^n$ with their root bubbles
$\{(M_\infty^j, g_\infty^j)\}_{j=1}^n$ and the associated centers
and scales $\{(p_i^j, \lambda_i^j)\}_{j=1}^n$. We define
$$
s^{m+1}_i(p) = r: \ \text{such that} \ \int_{B^i_r(p) \setminus
\bigcup_{j=1}^n B^i_{K^j\lambda_i^j}(p_i^j)} |Rm^i|^2 = \frac \delta
2 \tag 4.28
$$
and
$$
\lambda_i^{m+1} = s^{m+1}_i(p_i^{m+1}) = \inf_{M_i \setminus
\bigcup_{j=1}^n B^i_{K^j\lambda_i^j}(p_i^j)} s^{m+1}_i(p). \tag 4.29
$$
Then we repeat the procedure starting at the beginning of section 4.
Either $(p_i^{m+1}, \lambda_i^{m+1})$ is separable from all bubble
trees $\{T^j\}_{j=1}^n$, or it is a parent of bubble trees from
$\{T^j\}_{j\in J} \subset \{T^j\}_{j=1}^n$, which will be called the
new root bubble of this tree, according to (4.12). Clearly,
$$
m \leq \frac {2\Lambda}{\delta}.
$$
Thus the procedure has to stop in finitely many  steps. When the
procedure stops, we have $l$ number of separable bubble trees
$\{T^k\}_{k=1}^l$ with the root bubbles \newline $\{(M_\infty^k,
g_\infty^k)\}_{k=1}^l$ and their associated centers and scales
$\{(p_i^k, \lambda_i^k)\}_{k=1}^l$. Thus for a fixed number $\nu
>0$,
$$
\int_{B^i_{\nu}(p_i) \setminus \bigcup_{k=1}^l
B^i_{K^k\lambda_i^k}(p_i^k)} |Rm^i|^2 < \frac \delta 2. \tag 4.30
$$
After finitely many formation of exotic intermediate bubbles as in
the second step of the proof of Lemma 4.4, we finally arrive at
bubble trees $\{T^j\}$ whose associated centers are separated from
each other by non-zero distances in $M_i$. Thus the neck that
connecting each bubble tree $T^j$ to $M_i$ is given again by
applying the neck Theorem 3.1 to $B^i_{\sigma}(p_i^j)\setminus
B^i_{K^j\lambda_i^j}(p_i^j)$ for any small fixed positive number
$\sigma$, where $(M_\infty^j, g_\infty^j)$ is the root bubble for
$T^j$. This completes the construction of the bubble tree at any
point of curvature concentration in $M_i$.

One consequence of the bubble tree construction presents an
alternative argument for the following important estimates in [An-2]
[TV-2]: \proclaim{Corollary 4.6} Suppose that $(M_i, g_i)$ is a
sequence of Bach flat 4-manifolds $(M_i, g_i)$ with strictly
positive scalar curvature Yamabe metrics, i.e.
$$
Y(M_i, [g_i]) = \inf_{g\in [g_i]} \frac {\int_{M_i} R_g
dv_g}{(\int_{M_i}dv_g)^\frac 12} \geq Y_0 >0, \tag 4.31
$$
for some fixed number $Y_0$, vanishing first homology, and
$$
\int_{M_i} |Rm^i|^2dv_i \leq \Lambda. \tag 4.32
$$
Suppose that $X_i\subset M_i$ including some geodesic ball of radius
$r_0$ satisfies
$$
\int_{T_{\eta_0}(\partial X_i)} |Rm^i|^2 dv^i \leq \frac \delta 4,
\tag 4.33
$$
for some fixed positive numbers $r_0 > 4\eta_0 > 0$, $\delta =
\min\{\tau_0, \tau_1, \delta_0\},$ where $\tau_k$ are from Lemma 2.7
and $\delta_0$ is from Theorem 3.1. Assume there exists only one
bubble tree $T$ representing the curvature concentration in $X_i$.
Denote the center of the root bubble $(M_\infty, g_\infty)$ for $T$
by $p_i$. Then there exist a small number $\sigma_0$ such that the
intrinsic diameter of $S_\sigma(p_i)$ in $X_i \setminus
B^i_\sigma(p_i)$ is bounded by $C_0 \sigma$ for any $\sigma \leq
\sigma_0$ and some fixed $C_0$.
\endproclaim

\demo{Proof} Let $(M_\infty, g_\infty)$ be the root bubble for the
bubble tree $T$ representing the only curvature concentration in
$X_i$ and $(p_i, \lambda_i)$ be the associated center and scale.

Let $K$ be so large that
$$
\int_{M_\infty\setminus B^\infty_{K}(p_\infty)} |Rm^\infty|^2
dv^\infty \leq \frac \delta 4
$$
and
$$
\Cal H^3(S^\infty_r(p_\infty)) \leq 2v_4 r^3
$$
for $S^\infty_r(p_\infty) \subset M_\infty$, $r\geq K$ and $v_2 =
4v_4$, where $v_2$ is the same as in Theorem 3.1. Then let $\mu_0$
is so small that
$$
\int_{B^i_{\mu_0}(p_i)\setminus B^i_{K\lambda_i}(p_i)} |Rm^i|^2 dv^i
\leq \frac \delta 2
$$
and the neck Theorem 3.1 is applicable to the annulus
$B^i_{\mu_0}(p_i)\setminus B^i_{K\lambda_i}(p_i)$. Thus take
$\sigma_0 = \delta_0^{-\frac 14}\mu_0$. The intrinsic diameter of
$S_{\sigma}(p_i)$ for $\sigma \leq \sigma_0$ is then bounded because
of the upper bound of volume by Theorem 3.1  and the lower bound of
volume given by the Sobolev inequality (2.2).
\enddemo

\vskip 0.1in\noindent{\bf 5. Proof of main theorems}\vskip 0.1in

In this section we will apply the bubble tree construction in the
previous section and the recent compactness results in [An-2] [TV-2]
to prove our main theorems.

We first recall a recent result in [An-2], [TV-2]:

\proclaim{Theorem 5.1}([An-2] [TV-2]) \quad Suppose that $(M_i,
g_i)$ is a sequence of Bach flat, strictly positive Yamabe constants
with uniform lower bound, $b_1(M_i) = 0$, uniform $L^2$-Riemannian
curvature bound (4.4), and normalized volume
$$
\text{vol}(M_i, g_i) = \frac {8\pi^2}3. \tag 5.1
$$
Then a subsequence of $(M_i, g_i)$  converges to a Bach flat
manifold $(M_\infty, g_\infty)$ with finitely many isolated
irreducible orbifold points in the Gromov-Hausdorff topology.
Moreover, away from a finite set of points of curvature
concentration, the convergence is in $C^\infty$.
\endproclaim

We now begin the proof of Theorem A in section 1.

\demo{Proof} We will first establish the result in the case that $M$
is diffeomorphic to $S^4$, whose Euler number is $2$. In this case,
by the Gauss-Bonnet-Chern formula (1.2), we see that under the
assumption (1.9), we have
$$
\int_M |W|^2 dv \leq 16\pi^2 \epsilon.
$$
Choose $g$ to be the Yamabe metric on $M$, apply Kato's inequality
and Sobolev inequality (2.2), then there exists some constant $C$
depending on the Yamabe constant, such that
$$
(\int_M |E|^4 dv)^\frac 12 \leq C \int_M (6|\nabla E|^2  +
R|E|^2)dv. \tag 5.2
$$
Thus applying (2.6) in Lemma 2.3, we have
$$
\aligned (\int_M |E|^4 dv)^\frac 12 & \leq C (\int_M |W||E|^2
dv + \int_M |E|^3 dv) \\
& \leq C ((\int_M|W|^2dv)^\frac 12 + (\int_M |E|^2dv)^\frac
12)(\int_M |E|^4 dv)^\frac 12. \endaligned  \tag 5.3
$$
Thus, there is a positive number $\epsilon_0$ such that $E=0$ if
$\epsilon \leq \epsilon_0$ in (1.9). When $E=0$, we may use a
similar argument involving Kato's and Sobolev inequalities  as above
and apply (2.7) in Lemma 2.3 to conclude $W =0$ provided that
$\epsilon_0$ is sufficiently small. Thus $(M,g) $ is isometric to
$(S^4, g_c).$

We will now prove that $M$ has to be diffeomorphic to $S^4$. To see
this, we first first apply Theorem 5.1 above, together with the
bubble tree construction in section 4 to conclude:

\proclaim{Lemma 5.2}
$$
\text{vol}(M_\infty, g_\infty) = \frac {8\pi^2}3. \tag 5.4
$$
\endproclaim

We will prove $M$ is diffeomorphic to $S^4$ by a contradiction
argument, Assume the contrary,  then there is a sequence of Bach
flat manifold $(M_i, g_i)$ satisfying all assumptions in Theorem A;
and hence the assumptions of Theorem 5.1 by Corollary 2.8, with
$$
\int_{M_i} (\sigma_2 dv)[g_i] = 16\pi^2 (1 - \epsilon_i) \to
16\pi^2, \tag 5.5
$$
and $M_i$ is not diffeomorphic to $S^4$. Apply Lemma 2.6 to (5.5),
we thus get
$$
\int_{M_i} (|E|^2 dv)[g_i] = 32\pi^2 \epsilon_i \to 0. \tag 5.6
$$
Apply Theorem 5.1, we then conclude that some subsequence of $(M_i,
g_i)$ converges in the Gromov-Hausdorff topology to a limiting space
$(M_{\infty}, g_{\infty})$ which is an Einstein orbifold satisfying
$Ric_{\infty} = 3 g_{\infty} $. Apply Lemma 5.2 and the Bishop
volume comparison result, we then conclude that $(M_\infty,
g_\infty)$ cannot contain any orbifold point, hence it is a smooth
Einstein manifold; and thus it is isometric to the standard
4-sphere.

We now trace through the bubble tree construction in section 4.
First we note that at each point of curvature concentration, it is
the vertex of a bubble tree. Apply (5.6), we conclude that each
bubble is a Ricci flat bubble.

Since the limit space is smooth, the neck which connects the limit
space and the root bubble has to be a cylinder $[a,b]\times S^3$.
Hence the root bubble is Ricci-flat asymptotically Euclidean instead
of just being an ALE space. We now follow the proof of the
well-known result any Ricci flat asymptotically Euclidean orbifold
has to be the Euclidean space ( c.f. the proof of Theorem 3.5 in
[An-1]). Hence each neck in the bubble tree has to be a cylinder
$[a,b]\times S^3$. This in turn implies that the leaf bubbles are
Euclidean space, this contradicts the fact that the leaf bubbles
carry curvatures. We have thus finished the proof of Theorem A.
\enddemo

\proclaim{Theorem 5.3} The collection $\Cal A$ of Riemannian
4-manifolds $(M^4, g)$ satisfying the following

1) Bach flat,

2)
$$
Y(M, [g]) = \inf_{g\in [g]} \frac {\int_M R_g dv_g } {(\int_M dv_g
)^\frac 12} \geq Y_0, \tag 5.7
$$
for a fixed $Y_0 >0 $,

3)
$$
b_1(M) =  0, \tag 5.8
$$

4)
$$
\int_M (|W|^2 dv)[g] \leq \Lambda_0
$$
for some fixed positive number $\Lambda_0$. Then there are at most
finitely many diffeomorphism types in the family $\Cal A$.
\endproclaim

We will first assume Theorem 5.3 and establish Theorem B in section
1 as a consequence.

\demo{Proof of Theorem B} It is easy to see that (1.10) implies
(5.7) for manifolds with positive Yamabe constant. Hence (5.8) is a
consequence of the vanishing theorem of Gursky [G] (see Theorem 2.6
for the statement of the theorem). This finishes the proof of
Theorem B in section 1.
\enddemo

\demo{Proof of Theorem 5.4}

Assume otherwise, there is an infinite sequence of manifolds $(M_i,
g_i)$ from the collection with pairwise distinct diffeomorphism
types. We will show that, at least for a subsequence, manifolds
$(M_i, g_i)$ will be diffeomorphic to each other for $i$ sufficient
large.
 Without loss of generality, we
may assume that each $(M_i, g_i)$ satisfies all assumptions in
Theorem 5.1. So some subsequence of $(M_i, g_i)$ converges to a
limit space $(M_\infty, g_\infty)$ in the Gromov-Hausdorff topology
and in $C^\infty$ away from a finite set of points of curvature
concentration. In addition, at each point $p_\infty^k$ of curvature
concentration of $M_\infty$, as shown in the previous section, there
is a bubble tree $T^k$ forming for possibly some
 subsequence of $(M_i,
g_i)$. Since the set of points of curvature concentration is finite
and the bubble tree $T^k$ at each point $p_\infty^k$ is a finite
tree, we see that all the constants $K$'s (see (4.7) and (4.15) ) in
the bubble tree construction  in section 4 is a fixed finite set for
the particular subsequence we are considering. To show that $M_i$
are all diffeomorphic to each other, at least, for $i$ sufficiently
large, we consider each of the following three different type of
regions:

\vskip 0.05in {\bf Body region}: Let $\{(p_i^k, \lambda_i^k)\}_{k\in
N}$ be all the root bubbles of all the bubble trees. Then
$$
( M_i \setminus \bigcup_{k\in N} B^i_{\frac 14 \delta_0^\frac 14
}(p_i^k), g_i)
$$
tends to
$$
(M_\infty \setminus \bigcup_{k\in N} B^\infty_{\frac 14
\delta_0^\frac 14}(p_\infty^k), g_{\infty})
$$

\vskip 0.05in {\bf Neck region}: \quad Let $(p_i^{j_2},
\lambda_i^{j_2})$ be the parent for $(p_i^{j_1}, \lambda_i^{j_1})$.
Then
$$
B^i_{\frac 12 \delta_0^{\frac 14} \lambda^{j_2}_i}(p_i^{j_1})
\setminus B^i_{2 \delta_0^{-\frac 14}
K^{j_1}\lambda_i^{j_1}}(p_i^{j_1})
$$
tends to an annulus on an cone $C(S^3/\Gamma)$ for some finite group
$\Gamma\subset O(4)$ with a uniform bound on the order of the group
$\Gamma$ according to Theorem 3.1.

\vskip 0.05in {\bf Bubble region}: \quad Let $\{(p_i^j,
\lambda_i^j)\}_{j\in J}$ be all the children of a bubble $(p_i^k,
\lambda_i^k)$ in the bubble tree. Then
$$
B^i_{4 \delta_0^{-\frac 14}K^k\lambda_i^k}(p_i^k) \setminus
\bigcup_{j\in J} B^i_{\frac 14 \delta_0^{\frac
14}\lambda_i^k}(p_i^j)
$$
tends to
$$
B^\infty_{4\delta_0^{-\frac 14}K^k}(p_\infty^k)\setminus
\bigcup_{j\in J} B^\infty_{\frac 14 \delta_0^{\frac 14}}(p_\infty^j)
\subset M_\infty^k \setminus \{p_\infty^j\}_{j\in J}.
$$
Notice that the overlap regions of any two of the three different
type regions above are also well controlled and $M_i$ is covered by
those regions of the above three types. So $M_i$'s are all
diffeomorphic to each other for sufficiently large $i$ in the
subsequence. This contradicts with our assumption that $M_i$ are
pairwise not diffeomorphic to another. We have thus finished the
proof of Theorem 5.3.
\enddemo

\vskip 0.1in\noindent{\bf 6. Application to the analysis of the
$\sigma _2$ equation.}\vskip 0.1in

In this section, we study the  $\sigma _2$ equation for conformal
structures in the class $\Cal A$. Given a conformal structure $[g]
\in {\Cal A}$, represented by a metric $g_0$, a conformal metric
$g=e^{2w}g_0$ has its Schouten tensor $A_g$ defined as
$$
A_g = \frac 12 (Ric_g - \frac 16 R_g g).
$$
Under conformal change of metrics, we have
$$
A_g= -\nabla_0^2 w + dw {\otimes} dw - \frac{|\nabla_0 w|^2}{2}g_0 +
A_{g_0}.\tag 6.1
$$

We recall in the following two propositions the local estimate
developed in [GW] (see also the recent simplified proof of the
estimates in [Ch]) and the classification of entire solutions
[CGY-3]. A solution $g=e^{2w} g_0$ to the equation $ \sigma_2 (g)
=\sigma _2(A_g)=1$ is said to be admissible if the scalar curvature
$R_g$ is positive.

\proclaim{Proposition 6.1} Let $w \in C^3$ be an admissible solution
of the equation $\sigma _2(A_g)=1$ in $B_r$. There is a constant $C$
depending on $||g_0||_{C^4(B_r)}$ such that
$$
\sup_{B_{r/2}} ( |\nabla_0 ^2 w| + |\nabla_0 w|^2 ) \leq
 C (1+ \sup_{B_r} e^{2w}).  \tag 6.2
$$
\endproclaim

\proclaim{Proposition 6.2} An entire admissible solution
$g=e^{2w}|dx|^2$ to the $\sigma_2 (A_g) =1$ equation on $R^4$ is the
pull back of the spherical metric under the stereographic
projection.
\endproclaim

\noindent We will also use the following result of Viaclovsky ([V],
Proposition 3)

\proclaim{Proposition 6.3} Suppose $g_1=e^{2w_1} g_0$ and
$g_2=e^{2w_2} g_0$ are two metrics with $ \sigma_2 (g_1) \geq c_1 $
and  $ \sigma_2 (g_2) \geq c_2 $ for some positive numbers $c_1$ and
$c_2$ and with $R_{g_1}$ and $R_{g_2}$ both positive, then the
metric $g_t = e^{2 w_t} g_0$ where $ e^{-w_t} \equiv (1-t) e^{-w_1}
+ t e^{- w_2 } $ has $R_{g_t}$ positive and $ \sigma_2 (g_1) \geq
c_t $ for some positive $c_t$ depends only on $c_1$, $c_2$ for all $
0\leq t \leq 1$.
\endproclaim

\noindent {\bf Proof of Theorem C} \quad We may assume the conformal
structures under consideration are distinct from the standard
4-sphere, and hence in view of Theorem A, there is some $\epsilon_0
>0  $ that the following holds:
$$
\int \sigma _2(A_g)dV_g \leq 16 \pi ^2 (1 - \epsilon_0) . \tag 6.3
$$

We will establish the theorem by a proof of contradiction. Suppose
there is a sequence of metrics $g_i$ with $ \sigma_2 (g_i) = 1$ in
the family $\Cal A$ and with the diameter of $g_i$ tends to
infinity; we apply the argument in this paper to form a bubble-tree
of a subsequence of the underlying Yamabe metrics $(g_i)_Y$ of
$g_i$. We will then prove that the diameter of the metric $g=g_i$ in
each of the "body" and "neck" region of this bubble tree is bounded
and get a contradiction. We start with a simplest case when the
underlying Yamabe metrics of the sequence $\{g_i\}$ is compact in
the $C^{\infty}$ topology.

\vskip .1in

\noindent \proclaim{Lemma 6.4} Suppose  $g = e^{2w} g_Y $ is a
family of admissible metrics with $ \sigma_2 (g) = 1$ satisfying
both (6.3) and
$$
\int \sigma _2(A_g)dV_g \geq a_0 >0  \tag 6.4
$$
with the underlying Yamabe metrics $g_Y$ compact in the $C^{\infty}$
topology. Then there is a uniform upper and lower bound of the
conformal factor $w$, and the Ricci curvature of the metric is
bounded.
\endproclaim
\vskip .1in

\demo{Proof} First we claim that $w$ is bounded from above by the
following simple blowup argument. Thus we assume that there exists
diffeomorphisms $\Psi _i :(X,\bar g) \rightarrow (X_i, (g_i)_Y)$ so
that the sequence of pull back metrics $\Psi _i ^* (g_i)_Y$
converges in the $C^{\infty}$ topology to the limit metric $ \bar
g$. If the sequence $w_i$ were not bounded from above, there is a
sequence of points $p_i \in X_i$ such that $e^{w_i(p_i)}=max \;
e^{w_i}=\lambda _i$ tends to $\infty$. The diffeomorphisms
$\Psi_i^{-1}$ maps this sequence to a sequence $x_i \in X$ with a
convergent subsequence (still denoted by $\{ x_i \}$ ) to a point
$x_0 \in X$. Let $B$ be a ball of radius $r_0$ in a geodesic normal
(with respect to the limit metric $\bar g$) coordinate system $y$
whose origin correspond to the limit point $x_0$. Let $y_i$ denote
the coordinates of $\Psi _i^{-1} (x_i)$ and $T_i(y)=\lambda _i^{-1}y
+ y_i$ be a family of dilations, and consider the sequence of
metrics $h_i=T_i^*\Psi _i ^* (g_i)_Y$. The metrics $h_i$ are
isometric to $(g_i)_Y$ but defined on balls of radius $\lambda _i
r_0$ in $y$ space. Since the metrics $\Psi _i ^* (g_i)_Y$ converges
to $\bar g$, it follows that $h_i$ converges in $C^{\infty}$
uniformly on compact subsets in $y$-space to the flat metric
$|dy|^2$. The conformal metrics $T_i^*\Psi _i ^* {g_i}$ are
isometric to $ {g_i}$ can be written as $\lambda _i^{-2} e^{2w_i
\circ \Psi _i} h_i=v_i h_i$, where $v_i$ is a bounded function on
its domain of definition which includes the ball $|y|< \lambda _i
r_0$. Proposition 6.1 then asserts that euclidean $y$-gradient
$|\nabla v_i|$ is uniformly bounded on compact subsets, and hence
the functions $v_i$ converges uniformly on compact subsets to a
function $v$ on the $y$-space where the metric ${ h}=v^2|dy|^2$ is
an entire solution of the equation $\sigma _2(A_{h})=1.$ Proposition
6.2 asserts such solutions are the standard spherical metric. This
means however that the original metrics ${ g_i}$ must satisfy
$\limsup \int \sigma _2(A_{ g_i}) dV_{ g_i} \geq 16 \pi ^2.$ This is
a contradiction to (6.3). Thus we have shown that the conformal
factor $e^w$ and hence $w$ is bounded from above. Denote $ \bar w \,
=\,  max \,\,  w$, we now observe as $A_g$ satisfying condition
(6.4), we have
$$
 0 < a_0 \leq \int \sigma_2 (A_g) dv_g \,= \,  \int dv_g = \int e^{4w} dv_{g_Y}
\leq e^{4 \bar w} vol (g_Y), \tag 6.5
$$
Thus $\bar w $ is also bounded from below; from this we can apply
the local gradient estimate in Proposition 6.1 to conclude that that
$w$ is bounded both from above and below, applying Proposition 6.1,
we then conclude that the Ricci curvature of $g$ is uniformly
bounded. We have finished the proof of the lemma.
\enddemo
\vskip .2in

We now consider the general situation when the family of conformal
structures namely the family of the Yamabe metrics $g_Y$
corresponding to the $\sigma_2$ metrics may degenerate. According to
the proof of Theorem B, there is at most a bounded number of neck
regions and body region in any degenerating family of conformal
structures. To establish Theorem C, it suffices to check that in any
neck region of an ALE space which is the blown up limit of the
corresponding Yamabe metrics in the family ${\Cal A}$, the
restriction of the conformal metric satisfying the equation
$\sigma_2 =1$ has a bounded diameter. At this point it is necessary
to make precise the terms body region and neck region of the bubble
tree. The bubble tree refers to the pattern of degeneration of a
sequence of Yamabe metrics $(M_i, h_i)$ in the family ${\Cal A}$.
Thus there is a finite number of disjoint regions $\Omega _{j,i}, j
\in \{ 1,2,.,.,K \} $ in each $(M_i, h_i)$ over which the rescaled
(with possibly different scaling for each region) metrics $h'_i$
converges uniformly smoothly to a scalar flat, Bach-flat metric
which we denote $(\Omega_j, \, h'_j)$ and will call the $j-th$ body.
There is also a finite number of regions $T_{k,i}$ which overlaps in
$M_i$ with a pair of the previously labeled regions $\Omega _{j,i}$
and $\Omega _{j',i}$ such that $(T_{k,i}, h'_i)$ converges uniformly
on compact and smoothly to a subset of the end of the ALE space, and
such regions will be called the neck region. We also refer the
readers to the more detailed discussion of the bubble tree
construction and the definitions of the "body" and "neck" regions in
section 5 of the paper.

We remark that over any "body" region of the bubble tree, the
underlying rescaled Yamabe metrics in the family ${\Cal A}$ is
compact in the $C^{\infty}$ topology. Thus we may apply the argument
in the proof of the Lemma above to see that $w$ is bounded above for
each body in the family of metrics in $\Cal A$, hence the diameter
of the metrics bounded above in each of the body region of the fixed
bubble tree.

Our second remark is that over a body region of a bubble tree, if
the $\sigma_2$ "mass" defined to be the integral of $\sigma_2 (A_g)$
over the body is bounded from below; that is if condition (6.5)
holds for some constant $c_0$ replacing $ a_0$, then Lemma 6.4 also
applies to metrics in the family in this region. We shall call such
a region with the metric $g$ of bounded geometry.

\vskip .1in Since the total $\sigma_2$ mass of metrics in $\Cal A$
is  greater than $a_0$, and the total number K of body regions and
neck regions is finite, we have two cases: \vskip .1in Case I: There
is some body region where the mass over there is greater than $
\frac {1}{K} a_0$. \vskip .1in Case II: Case I does not happen,
while there is a neck region where the mass is greater than $\frac
{1}{K} a_0$. \vskip .1in In the following, we will first Analyze the
situation in case II. In this case, as we are studying the metric
$g$ over a region which asymptotically looks like $(S^3/\Gamma)
\times R$, where $ \Gamma$  is a finite group with the order $
\|\Gamma \|$ uniformly bounded for metrics in the family $\Cal A$,
we will use the cylindrical coordinate $x=e^t \nu ;\; t\in R;\; \nu
\in (S^3/\Gamma)$ and choose the background metric $g_0$ to be the
standard cylindrical metric $ g_0 = dt^2 + d \nu^2$. Denote $ g =
e^{2w} g_0$; we may apply argument similar as before and prove that
$w$ is uniformly bounded above for $g \in \Cal A$. In the following
we will show that in case II, we also have  max w is also bounded
from below by a constant depending only on the $\sigma_2$ mass over
the region and K. Applying Proposition 6.1, we then conclude that
there is a region (around the point where max w occurs) where the
metric $g$ has locally bounded geometry.

To prove such a result for a metric $g = e^{2w} g_0$ in $\Cal A$
over a neck region with $\sigma_2(g) =1$, we first observe that we
may replace $ w = w(t, \nu)$ by $ \bar w (t) = \frac
{1}{|S^3|/|\Gamma| } \int_{S^3/\Gamma} w(t, \nu) d \nu $ and apply
Proposition 6.3 to conclude that the metric $ \bar g = e^{ 2 \bar w}
g_0$ still satisfies that $\sigma_2 (\bar g) \geq c_3 > 0$ and $
R_{\bar g}$ positive, and where $c$ is a constant. Notice that in
this case, we have $ |w (t, \nu) - {\bar w} (t)| $ being uniformly
bounded, thus
$$
\int e^{ 4 \bar w (t)} dt \geq C \int e^{ 4  w (t, \nu)} dt d\nu
\geq C(K) \sigma_0.
$$
Thus it suffices to establish the following Lemma. The idea of the
proof of the Lemma follows from the analysis of O.D.E. solutions of
$\sigma_k (g) = constant$ solutions on annulus regions in the
earlier work of [CHY].

\proclaim{Lemma 6.5} Suppose  $g = e^{2w} g_0 $ where $w = w(t) $ is
a metric $\sigma_2 (\bar g) \geq c_3 > 0$ and $ R_{\bar g}$ positive
defined on a cylinder $(S^3/\Gamma) \times [ t_0, T ]$ with
$$
\int_{t_0} ^T e^{4w (t)} dt = a > 0  \tag 6.6
$$
then
$$
max w \geq  \frac 12 log a + c_4
$$ for some universal constant $c_4$
depending only on $c_3$.
\endproclaim

\demo{Proof} In the cylindrical coordinate, we have
$$
\sigma_2 (A_g)(t)  = - \frac 23 w'' (1- (w')^2) e^{-4w} (t) \geq c_3
\tag 6.7
$$
while $R_g (t) = 6 ( 1 - (w')^2 - w'')(t) e^{ - 2 w(t)} >0$. Thus we
have $ (1- (w')^2)>0 $ while $ w'' < 0$; thus $ w'$ is a decreasing
function. Suppose $ \, max \, w \,  = \,  w (t_M)$  happens for some
point $t_M$ in ($t_0$, T), then on one of the interval ($t_0$,
$t_M$), or ($t_M$, T) the $\sigma_2$ mass is greater than $\frac 12
a $. Choose the interval this happens, say $(t_M, T)$. Otherwise, we
may assume by reversing the variable $t$ to $-t$ if necessary, that
$ \, max \,  w \, = \, w(t_M)$. We now fix a number $ b = \frac 12
\int_{t_M}^T \sigma_2 (g) dt $, and choose N so that
$$
\int_{t_M}^{t_M + N} e^{4w(t)}  dt = b.
$$
Now for $t> T_M + N$, we have from (6.7)
$$
- w''(t) \geq c e^{ 4 w(t)},
$$
for $c =\frac 32 c_3$. Thus
$$
-w'(t) \geq c \int_{t_M}^t e^{ 4 w(s)}ds  \geq cb,
$$
hence a direct integration yields the bound
$$
b = \int_{t_M+N}^T e^{ 4 w(t) } dt \leq \frac { e^{ 4 w (t_M) }}
{4bc}. \tag 6.8
$$
The assertion of the lemma is a direct consequence of (6.8).
\enddemo
\vskip .2in

We remark that as a consequence of Lemma 6.5, when condition (6.6)
is satisfied for some $ a \geq \frac {a_0}{K} $, then the region $B
\, = \, \{ (t, \nu): | t - t_M| \leq 1, \nu \in  S^3/\Gamma \} $ has
g-volume bounded from below by a positive constant depending only
the data in the family $\Cal A$; and in fact is a region where
$Ric_g$ bounded from above and below by positive constant.

Thus we have arrived at the conclusion that in any degenerating
sequence of $\sigma _2=1$ conformal metrics, either in case I or in
case II discussed above, there is a region of bounded geometry
namely the body region in case I, the region in the cylinder which
has mass in case II which we denote by $B$. We will now use the
volume comparison argument to show that the diameter of each neck is
necessarily bounded. We proceed to argue by contradiction. Assume
there is one neck region over which the sequence of the $\sigma
_2=1$ metrics have unbounded diameter. Let $d$ denote the distance
function in the $g$ metric to a point $p$ to be specified later.
Since $\int _{\Omega} \sigma_2 dV$ is given by the volume of the
neck region $\Omega$ and is bounded above by $16\pi ^2$, it follows
that given any constant $\epsilon_0 > 0$, there is a length $L$ so
that the geodesic annulus $A_{r,  r+L}$ contains a geodesic sphere
$S_{\rho}$ whose volume is bounded above by $\epsilon _0$. Now as a
$\sigma_2 =1$ metric has nonnegative Ricci curvature, according to
the volume comparison result of Bishop-Gromov, we have
$$
\frac {Vol(S_{\rho + s})}{( \rho + s)^3} \leq \frac{Vol(S_\rho
)}{{\rho}^3}, \tag 6.9
$$
for each $ s>0$. It is necessary at this point to remark that in the
inequality above, the geodesic spheres $S_{\rho + s}$ and $S_{\rho}$
may be replaced by a component of the geodesic sphere denoted by
$S'_{\rho + s}$ and $S'_{\rho}$ provided the length minimizing
geodesic joining $p$ to points in the given component of $S'_{\rho +
s}$ passes through the corresponding component in $S'_{\rho}$. Thus
if $s  \leq 2L$ and $\rho$ is much larger than L, we find
$Vol(S'_{\rho + s}) \leq 2 \epsilon _0$ since each point $q$ in the
neck is at $g$ distance no more than $L$ from such a geodesic sphere
$S'_{\rho} $ if we choose to measure distance from a point $p$
located at the far end of the neck from $q$. Thus we have shown that
each component of the geodesic sphere $S'_{\rho}$ has uniformly
small volume as long as the neck has sufficiently large diameter.
Now suppose $p$ is chosen to lie in the far end of the neck (away
from the region of bounded geometry, Denote $D= d(p, B)$, then the
component of spheres $S'_{D+\lambda}$ as $\lambda$ increases from
zero to diameter $d$ of $B$ will sweep out the region $B$. Using the
Fubini theorem, In particular for some choice of $\lambda$ the set
$S'_{ D+ \lambda }$ will have area bounded from below by $volume
(B)/d$. This is in contradiction to the uniform smallness of the
volume of all such component of spheres. Thus we have finished the
proof of Theorem C.
\enddemo

\vskip 0.1in \noindent {\bf References}:

\roster
\item"{[An-1]}" M. Anderson; Ricci curvature bounds and Einstein
manifolds on compact manifolds, JAMS. vol. 2, no. 3, (1989),
455-490.

\vskip 0.1in
\item"{[An-2]}" M. Anderson; orbifold compactness for spaces of
Riemannian metrics and applications, preprint, arXiv:math.DG/0312111

\vskip 0.1in
\item"{[AC]}" M. Anderson and J. Cheeger; Diffeomorphism
finiteness for manifolds with Ricci curvature and $L^{n/2}$-norm of
curvature bounded, Geom. and Funct. Anal. vol. 1, no. 3, (1991),
 231-252.

\vskip 0.1in
\item"{[Au]}" T. Aubin; Equations
differentielles non lineaires et probleme de Yamabe concernant la
courbure scalaire, J. Math. Pures Appl. 55 (1976), 269--296

\vskip 0.1in
\item"{[Ba]}" S. Bando; Bubbling out of Einstein manifolds,
Tohoku Math.J.,II.Ser. 42, no.2, (1990), 205-216.

\vskip 0.1in

\item"{[BC]}" H. Brezis and J.M. Coron; Convergence of solutions of
H-systems or how to blow bubbles, Archive Rat. Mech. Anal. 89
(1985), 21-56.

\vskip0.1in
\item"{[BN]}" H. Bray and A. Neves; Classification of
prime 3-manifolds with Yamabe invariant larger than $RP^3$, reprint
2003, to appear in Annals of Math.

\vskip 0.1in
\item"{[CGY-1]}" S.-Y. A. Chang, M. Gursky and P. Yang;  An equation of
Monge-Ampere type in conformal geometry and 4-manifolds of positive
Ricci curvature, Annals of Math. 155 (2002), 709-787.

\vskip 0.1in
\item"{[CGY-2]}" S.-Y. A. Chang, M. Gursky and P. Yang; A conformally invariant
Sphere theorem in four dimension, Publications Math. Inst. Hautes
Etudes Sci. 98 (2003), 105-143.

\vskip 0.1in
\item"{[CGY-3]}" S.-Y. A. Chang, M. Gursky and
P. Yang; An apriori estimate for a fully nonlinear equation on
4-manifolds", Jour. D'Anal. Math., 87 (2002), 151-186.

\vskip 0.1in
\item"{[CHY]}" Sun-Yung A. Chang, Zheng-Chao Han and Paul Yang,
``Classification of singular radial solutions  to the $\sigma_k$
Yamabe equation on annular domains'',  JDE, 216 (2005), pp 482-501.

\vskip .1in
\item"{[C]}" S. Chen, "Local estimates for some fully nonlinear
elliptic equations', IMRN 2005, no. 55, pp 3403-3425.

\vskip .1in
\item"{[De]}" A. Derdzinski; Self-dual K\"hler manifolds and Einstein
manifolds of dimension four, Composition Math. 49 (1983), 405-433.

\vskip0.1in
\item"{[G]}" M. Gursky;
The Weyl functional, deRhan cohomology, and Kahler-Einstein metrics,
Annals of Math., 148 (1998), 315-337.

\vskip 0.1in
\item"{[GW]}" P. Guan and X.-J. Wang, "Local estimates for a class of
fully nonlinear equations arising from conformal geometry", IMRN
2003, no. 26, pp 1413-1432.

\vskip0.1in
\item"{[He]}" E. Hebey; ``Sobolev spaces on Riemannian
manifolds'', Lecture Notes in Mathematics, 1635 Springer-Verlag,
Berlin, 1996.

\vskip .1in
\item"{[Ma]}" C. Margerin; A sharp characterization of the smooth
4-sphere in curvature forms, CAG, 6, no. 1, (1998), 21-65.

\vskip 0.1in
\item"{[Q]}" J. Qing,  On Singularities of the Heat
Flow for Harmonic Maps from Surfaces into Spheres, Comm in Anal. and
Geo. vol. 3 (1995),  297--315.

\vskip 0.1in
\item"{[St]}" M. Struwe, Global compactness result for elliptic
boundary value problem involving limiting nonlinearities, Math. Z.
187 (1984), 511-517.

\vskip .1in
\item"{[S]}" R. Schoen; Conformal deformation of a Riemannian
metric to constant scalar curvature, J. Diff. Geom.  20 (1984),
479--495.

\vskip .1in
\item"{[TV-1]}" G. Tian and J. Viaclovsky, Bach flat asymptotically
ALE metrics, preprint, arXiv:math.DG/0310302.

\vskip 0.1in
\item"{[TV-2]}" G. Tian and J. Viaclovsky, Moduli space of critical
Riemannian metrics in dimension 4, preprint, arXiv:math.DG/0312318.

\vskip 0.1in
\item"{[V]}" J. Viaclovsky, Conformally invariant Monge Ampere equations:
global solutions, TAMS 352 (2000) 4371-4379.

\endroster
\enddocument